\documentclass[12pt]{amsart}
\usepackage{amsfonts,amssymb,amscd,amsmath,enumerate,verbatim,calc,graphicx}

\usepackage{thumbpdf} 

\usepackage[pagebackref]{hyperref} 

\theoremstyle{plain}
\newtheorem{Theorem}{Theorem}
\newtheorem{Lemma}[Theorem]{Lemma}
\newtheorem{Corollary}[Theorem]{Corollary} 
\newtheorem{Proposition}[Theorem]{Proposition}

\newtheorem*{Theorem-introduction}{{Theorem}}
\newtheorem*{Proposition-introduction}{{Proposition}} 

\theoremstyle{definition}

\newtheorem{ex}[Theorem]{Example}
\newtheorem{remark}[Theorem]{Remark}

\newcommand{\mlabel}[1]%
  {\mbox{}\marginpar{\raggedleft\hspace{0pt}{\rm\ttfamily#1}}\label{#1}}

\newcommand{\bR}{{\mathbb{R}}}
\newcommand{\bN}{{\mathbb{N}}}

\newtheorem{pr}{Problem}
\newcounter{hours}\newcounter{minutes}

\def\blfootnote{\xdef\@thefnmark{}\@footnotetext}  

\title[Intersection bodies]{An extension of a result by Lonke to  intersection bodies}
\author[M. Angeles Alfonseca]{M. Angeles Alfonseca}
\address{Department of Mathematics, North Dakota State University}
\email{maria.alfonseca@ndsu.edu}

\thanks{Work partially supported by the NDSU Advance FORWARD program sponsored by the National Science Foundation, HRD-0811239}

\begin{document}

\begin{abstract} 
In this paper we prove that intersection bodies cannot be direct sums using Fourier analytic techniques. This extends a result by Lonke. We also  prove a necessary regularity condition and a convexity condition for a body of revolution to be an intersection body of a star body. 

\end{abstract}
\maketitle

\section{Introduction}

Let $K,L$ be origin-symmetric star bodies in $\bR^n$. The body $K$ is called the {\it intersection body of $L$}, and denoted by  $K=I \! L$, if the radius of $K$ in each direction is equal to the $(n-1)$-dimensional volume of the central section of $L$ that is perpendicular to that direction. In other words, if $\rho_K(\xi)=\max\{a: a\xi \in K\}$ is the radial function of $K$, then for every $\xi \in S^{n-1}$, we have  $ \rho_K(\xi)={\mbox{Vo}}l_{n-1}(L\cap \xi^\perp)$.

The volume of the section of $L$ can be written using the spherical Radon transform $R$  (see \cite{G,K}):
\[
      \rho_K(\xi)=\frac{1}{n-1} \int_{S^{n-1}} \rho_L^{n-1} (\theta) \, d \theta= \frac{1}{n-1}\, R( \rho_L^{n-1})(\xi),           
      \qquad \forall \xi\in S^{n-1}.
\]
The more general class of {\it intersection bodies} is defined as follows. A star body $K$ is an intersection body if its radial function $\rho_K(\xi)=\max\{a: a\xi \in K\}$ is the spherical Radon transform of an even non-negative measure $\mu$, {\it i.e.}, for every continuous function $g$ on $S^{n-1}$,
\[
      \int_{S^{n-1}} \rho_K(\xi) g(\xi) \, d\xi = \int_{S^{n-1}} Rg(\xi) \mu(d\xi), 
              \qquad \forall \xi\in S^{n-1}.
\]
Intersection bodies were originally introduced by Lutwak \cite{Lu} in connection to the Busemann-Petty problem (see, for example, \cite{K}), and were instrumental in finding its complete solution \cite{GKS}.

A Fourier-analytic characterization of intersection bodies is due to Koldobsky \cite{K}. If $\|\xi\|_K=\rho_K(\xi)^{-1}$ denotes the norm associated to the  body $K$, then $K$ is an intersection body if and only if the Fourier transform of $\| \cdot \|_K^{-1}$ is a positive distribution, {\it i.e.} its action on any non-negative test function gives a non-negative result. 

The structure and the geometric properties of intersection bodies are hard to understand  for at least two reasons. First, if $n=2,3,4$, all origin-symmetric convex bodies in $R^n$ are intersection bodies, and thus their study is only meaningful in dimension $n \geq 5$. Secondly,  Zhang proved that no polytope is an intersection body of a star body if $n\geq 4$ \cite{Z}.

Intersection bodies are closely related to zonoids, a very symmetric class of convex bodies. Zonoids admit several different characterizations  \cite{B}. For example, they are Hausdorff limits of sums of segments. Full dimensional zonoids are also centered projection bodies ($K$ is the projection body of $L$ if  the width of $K$ in every direction $\xi$ equals twice the $(n-1)$-dimensional volume of the projection of $L$ orthogonal to $\xi$). Notice the parallel with the definition of an intersection body of a star body. In a certain sense, it may be expected that the geometric properties of intersection bodies are dual to the geometric properties of zonoids, since there is a certain duality between sections and projections (see Section 7.4, Note 4  in \cite{Sc}). 
However, this duality does not hold in all instances. For example,  if $Z$ is a zonoid, then $Z^*$ is an intersection body, but   the class of intersection bodies is wider than the class of duals of zonoids.  A better understanding of intersection bodies could be useful in the solution of the following open problem.
\begin{pr}\label{pr:zongr}
 Is it true that all zonoids whose polars are zonoids tend to the Euclidean  ball when  $n \to \infty$? In other words, if $d_n$ is  the Banach-Mazur distance, is it true that $c=\lim\limits_{n\to \infty}d_n(Z, B^2_n)=1$?
\end{pr}

\noindent This problem is related to an isometric analogue of a well-known theorem of Grothendieck  \cite{LP}, asserting that among infinite dimensional Banach spaces, the ones which are isomorphic to both a subspace of $L^1$ and a quotient-space of $L^{\infty}$ are isomorphic to a Hilbert space.  R. Schneider \cite{S} used high-order derivatives of the support function to construct examples in all dimensions of zonoids whose polars are also zonoids, that are not ellipsoids. It is known that all his zonoids tend to the Euclidean ball as $n$ goes to infinity. It is not known if there are zonoids whose polars are zonoids that do not converge to the Euclidean ball. 

Y. Lonke  \cite{L} showed that  if $K=A+B$ is a convex body in $R^n$, $n\ge 3$, such that $dim \mbox{(span A)}=n$, $1\le dim \mbox{(span B)} \le n-2$, then the polar of $K$ is not a zonoid. This shows, in particular, that a zonoid whose polar is also zonoid cannot have an $(n-2)$-dimensional face. In the case of a direct sum, he was able to get rid of the restriction on the dimension of the face, and proved the result for $1\le dim \mbox{ (span B)} \le n-1$. On the other hand, Lonke proved that the barrel $Z_n=B_n^2+B_{n-1}^2$ is a zonoid whose polar is also a zonoid in dimensions $n=3,4$. Thus, it is possible for a zonoid whose polar is a zonoid to have an $(n-1)$- dimensional face. The construction of such examples in higher dimensions would give a negative answer to Problem 1.

In Section 2  we present a necessary condition for a body $K$ to be an  intersection body. It is an extension of Lonke's results, proved using the Fourier-analytic characterization and  Koldobsky's Second Derivative Test \cite{K5}.  We show that {\it if  $K\subset R^n$, $n \geq 5$ is a convex body that can be written as a {\it direct} sum $K=A+B$, where dim$(A) \geq 1$ and dim$(B)\geq 4$, then $K$ is  not an intersection body}. This was previously proved by Zhang, using the Radon transform (see Note 8.1 in \cite{G}). Note that such a result is not true for non-direct sums, as there exist intersection bodies with low dimensional faces. For example, if $L$ is  a body of revolution $L$ with a cylindrical part near the equator, then $I \! L$ has an $(n-1)$-dimensional symmetric face (\cite{G}, see the proof of Theorem 8.1.18).

In Section 3 we prove a regularity condition for a body of revolution to be the intersection body of a star body. We also find an equator-convexity condition that allows us to determine, given an intersection body of a star body of revolution in $\bR^n$, if the body in $\bR^{n+2}$ that has the same 2-dimensional radial function is still an intersection body in dimension $n+2$. The original motivation for the work in this section was to prove that Lonke's zonoid $Z_n$ is not an intersection body in dimension 5 and higher, thus explaining why its polar is not a zonoid in those dimensions. As it turned out, $Z_n$ is an intersection body in dimensions 5 and 6, but not in dimensions 8 and higher. 

Acknowledgements: The author would like to thank Dmitry Ryabogin and Artem Zvavitch for many conversations and ideas about this paper. The author also thanks the reviewer for useful suggestions to make the paper accurate and more readable.  

\section{Direct sums are not intersection bodies for $n \geq 7$}

Let $K$ be an origin-symmetric convex body in $\bR^n$, $n \geq 5$, that can be written as a direct sum $K=A \oplus B$, where $j=dim(A) \geq 1$ and $n-j=dim(B) \geq 4$. Note that, because of the direct sum, $A$ and $B$ inherit central symmetry from $K$. We shall show that $K$ is not an intersection body. Since a central section of an intersection body is also an intersection body \cite{GW}, it is enough to consider the case in which $j=1$, {\it i.e.} $A$ is a segment. We will write the points $z\in \bR^n$ in the form $z=(x,y)$, where $x \in \bR$ and $y \in \bR^{n-1}$. Assuming that the segment $A$ has length 2, the norm associated to $K$ can be written as
\[
  \|(x,y)\|_K=\max\left\{|x|,\|y\|_{B} \right\}.
\]

We shall use the second derivative test introduced in \cite{K5} to prove that $K$ is not an intersection body. Theorem 1 in \cite{K5} requires that the function $x \rightarrow \|(x,y)\|$ has continuous second derivative everywhere on $\bR$. However, in our calculations all the derivatives of the norm will be taken in the sense of distributions and we will not need this regularity. 

The following proof follows that of Lemma 1 in \cite{K5}. For every $m \in \bN$, we consider the functions $h_m(x)=\frac{m}{\sqrt{2 \pi}} e^{-x^2 m^2/2}$ and $u(y)=\frac{1}{(2 \pi)^{(n-1)/2}} e^{-\|y\|_2^2 /2}$.

\begin{Lemma}
 \label{cylinder}
Let $\|(x,y)\|$ be the norm defined by $K$. For every $\epsilon >0$ there exists $M \in \bN$ so that, for every $m>M$, $\, \langle\|(x,y)\|^{-1},u(y)\,h_m''(x)\rangle \geq -\epsilon$.

\end{Lemma}

\noindent{\bf Proof.} Let us define the sets 
\begin{eqnarray*}
U=\left\{
(x,y): \, |x| > \|y\|_B
\right\},\\ 
W=\left\{
(x,y): \, |x| < \|y\|_B
\right\}.
\end{eqnarray*}
Then
\[
  \langle\|(x,y)\|^{-1},u(y)h_m''(x)\rangle=
  \int_{U} u(y)\, h_m''(x) \, \|(x,y)\|^{-1} dx\,dy + 
  \int_{W} u(y) \, h_m''(x) \, \|(x,y)\|^{-1} dx\,dy
\]
If $(x,y) \in U$, then
\[
 \int_U u(y)\, h_m''(x) \, \|(x,y)\|^{-1} dx\,dy=
\]
\begin{equation}
 \label{integrala}
  \int_{\bR^{n-1}\setminus\{0\}} u(y) \int_{|x|>\|y\|_B} |x|^{-1}\, h_m''(x) dx \, dy=
\end{equation}
\[
  2 \int_{\bR^{n-1}\setminus\{0\}} u(y) \int_{\|y\|_B}^\infty \frac{h_m''(x)}{x} dx \, dy
\]
Integrating by parts twice, we obtain
\begin{equation}
 \label{3terms}
   2 \int_{\bR^{n-1}\setminus\{0\}} u(y)\left[
    -\frac{h_m'(\|y\|_B)}{\|y\|_B} 
    -\frac{h_m(\|y\|_B)}{\|y\|^2_B} 
    +2 \int_{\|y\|_B}^{\infty} \frac{h_m(x)}{x^3} dx
  \right]
\, dy.
\end{equation}
Since $u(y)$ and $h_m(x)$ are positive functions, and $h_m'(x)$ is negative, the first and third integrals in (\ref{3terms}) are positive. Thus, we only need to study the second integral,
\[
  -\int_{\bR^{n-1}\setminus\{0\}} u(y) \frac{h_m(\|y\|_B)}{\|y\|^2_B} \, dy
\]
Using the change of variables $z=my$, we can rewrite this integral as
\[
 -m^{4-n}\int_{\bR^{n-1}\setminus\{0\}} u\left(\frac{z}{m}\right)
 \frac{h_1(\|z\|_B)}{\|z\|^2_B} \, dz.
\]
Since the integral is well defined and $n \geq 5$, this term converges to 0 as $m$ goes to infinity, as we want.

Now we turn to the second case. If $(x,y) \in W$, then 
\begin{equation}
 \label{cylindrical}
 \int_{W} u(y)\, h_m''(x) \, \|(x,y)\|^{-1} dx\,dy=
\end{equation}
\[
  \int_{\bR^{n-1}\setminus\{0\}} u(y) \int_{|x|<\|y\|_B} \frac{h_m''(x)}{\|y\|_B} dx \, dy=
\]
\[
  2 \int_{\bR^{n-1}\setminus\{0\}} \frac{u(y)}{\|y\|_B} \int^{\|y\|_B}_0 h_m''(x)
  dx \, dy =
\]
\[
     2 \int_{\bR^{n-1}\setminus\{0\}} \frac{u(y)}{\|y\|_B}
    \left[h_m'(\|y\|_B)-h_m'(0)\right] \, dy =
\]
\[
   2 \int_{\bR^{n-1}\setminus\{0\}} \frac{u(y)}{\|y\|_B}
    h_m'(\|y\|_B) \, dy = -2m^2 \int_{\bR^{n-1}\setminus\{0\}} u(y) h_m(\|y\|_B) \, dy
\]
The change of variables $z=my$ gives
\[
   -2 \int_{\bR^{n-1}\setminus\{0\}} m^{4-n} u\left(\frac{z}{m}\right) h_1(\|z\|_B) \, dz
\]
which, as before, converges to $0$ when $m$ goes to infinity.

\qed

\bigskip

It follows (cf. \cite{K5}) that $K$ is not an 
intersection body. We write the proof for completeness. 

\begin{Theorem}

\label{direct}

Let $K$ be a direct sum of a segment and an $(n-1)$-dimensional body, $n\geq 5$. Then $K$ is not an intersection body. 

\end{Theorem}

\noindent{\bf Proof:} Suppose that $K$ is an intersection body. Then $\|(x,y)\|_K^{-1}$ is a positive definite distribution (by Theorem 4.1  of \cite{K}). This implies, by Lemma 2.24 and Corollary 2.26 in \cite{K}, that there exists a finite Borel measure $\mu_0$ on $S^{n-1}$ such that, for every even test function $\phi$, 
\begin{equation}
 \label{koldobsky}
    \int_{\bR^n} \|(x,y)\|_K^{-1} \phi(x,y) \, dx \,dy=  \int_{S^{n-1}} \left(   \int_0^\infty \widehat{\phi} (t \xi) \, dt \right)  d \mu_0(\xi).
\end{equation}
In our case, (\ref{koldobsky})  becomes
\begin{equation}
 \label{derivative}
    \int_{\bR^n}  \|(x,y)\|_K^{-1} \frac{\partial^2 \phi}{\partial x^2} 
   =- \int_{S^{n-1}} \xi_1^2 d \mu_0 (\xi) \, 
   \int_0^\infty  t^2 \widehat{\phi} (t \xi) \, dt.
\end{equation}
Let $\phi(x,y)=h_m(x)u(y)$, where $h_m(x)$ and $u(y)$ are the functions defined before Lemma 1. Then,
\[
    \widehat{\phi(\xi)}=e^{-\xi_1^2/2m^2}e^{-(\xi_2^2+\cdots+\xi_n^2)/2}.
\]
With this choice of $\phi$, (\ref{derivative}) becomes
\[
   -\epsilon \leq \langle (\|(x,y)\|_K^{-1}),h_m(x)'' \cdot u(y) \rangle
\]
\[ = -\sqrt{\pi/2} 
   \int_{S^{n-1}} \xi_1^2 \left(\frac{\xi_1^2}{m}+\xi_2+\ldots+\xi_n^2 \right)^{-3/2}
   \, d \mu_0 (\xi)
\]
\[
    \leq  -\sqrt{\pi/2} \int_{S^{n-1}} \xi_1^2 \, d \mu_0 (\xi) \leq 0.
\]
In the first inequality we have used Lemma \ref{cylinder}. Thus, the measure $\mu_0$ is supported on $S^{n-1}\cap \{\xi_1=0 \}$, which is a contradiction with the fact that $K$ is an $n$-dimensional body. 

\qed

\begin{remark}  Theorem \ref{direct} proves that any direct sum in $n \geq 7$ is not an intersection body, since one of its summands has at least dimension 4. In order to prove the same result for $n \geq 5$, there is only one case left to consider: a direct sum of a 2-dimensional and a 3-dimensional body in $\bR^5$.
\end{remark}

\begin{ex} The cylinder $B_{n-1}(0,1) \times [0,1]e_n$ is not an intersection body for any $n \geq 5$.
\end{ex}

\begin{ex} 
A slight variation of the proof of Lemma 1 allows us to prove a version of the Second Derivative Test for bodies of revolution in $\bR^n$  with an $(n-1)$-dimensional face. In general, such bodies  may or may not be intersection bodies (the cylinder in dimension 5 is not an intersection body, but the intersection body of the cylinder is). Our interest in knowing if this type of bodies are intersection bodies is related to Problem 1 mentioned in the Introduction. The construction in all dimensions of zonoids whose polar are zonoids, and that have lower dimensional faces, would provide a negative answer to Problem 1. Such bodies would necessarily be intersection bodies. Proposition \ref{revoll} below gives a necessary condition for a body of revolution with an $(n-1)$-dimensional face to be an intersection body.

Let $K \in \bR^n$ be a  body of revolution  with an $(n-1)$-dimensional face.
As before, we denote the points of $\bR^n$ as $(x,y)$, with $x\in \bR$ and $y \in \bR^{n-1}$. Assume that the axis of revolution of $K$ is the $x$-axis, and that the face perpendicular to it has radius 1 and is placed at height 1. Then the norm $\|\cdot\|_K$ can be written as
\[
    \|(x,y)\|_K=\left\{ 
 \begin{array}{cc}
    |x| & |x|>\|y\|_2 \\
    g(x,\|y\|_2) & |x|<\|y\|_2
 \end{array}
 \right.
\]
where $g(x,r)$ is positive, convex, homogeneous of degree 1 and even with respect to each variable. In particular, for every $r\neq 0$ fixed, $g(x,r)$ has a positive minimum at $x=0$. 
 \end{ex}

\begin{Proposition}
  \label{revoll}

Let $K$ be a body of revolution with a face and assume that: 
\begin{enumerate} 
  \item For every fixed $r\neq 0$, the function $x \rightarrow g(x,r)$ has continuous second derivative on $|x|<|r|$ and $\frac{\partial^2 g}{\partial x^2} (0,r)=0$. 
 \item $\lim_{x \rightarrow 0} \frac{\partial^2 g}{\partial x^2}(x,r)=0$ uniformly on $r$. 
\end{enumerate}
Then $K$ is not an intersection body.

\end{Proposition}

\noindent{\bf Proof:} 

With $u(y)$ and $h_m(x)$ defined as in Lemma 1, we consider the integral
\[
  \langle\|(x,y)\|_K^{-1},u(y)h_m''(x)\rangle=
\]
\[
  \int_{\bR^{n-1}\setminus \{0\}} u(y)\, \int_{|x|>\|y\|_2} 
   h_m''(x) \, |x|^{-1} dx\,dy + 
  \int_{\bR^{n-1}\setminus \{0\}} u(y)\, \int_{|x|<\|y\|_2} 
   h_m''(x) \, (g(x,\|y\|_2))^{-1} dx\,dy.
\]
The first of these two integrals is exactly the same as the integral in (\ref{integrala}). As we showed in (\ref{3terms}), it is well defined and converges to 0 as $m$ goes to infinity if $n \geq 5$. 

We integrate by parts the second integral and we use that $g$ is homogeneous of degree 1 and that $g(1,1)=1$ (from the definition of the norm of $K$), obtaining
\[
   2  \int_{\bR^{n-1}\setminus\{0\}} u(y)\left[ \frac{h_m'(\|y\|_2)}{\|y\|_2}
    + \int_0^{\|y\|_2}  h_m'(x) \frac{\partial g}{\partial x}(x,\|y\|_2) 
     \frac{1}{g(x,\|y\|_2)^2}  dx
  \right]
\, dy =I+I\!I.
\]
Making the change of variables $w=my$, $I$ is equal to
\[
 I=-\frac{2}{m^{n-4}}
 \int_{\bR^{n-1}\setminus\{0\}} u\left(\frac{w}{m}\right) h_1(\|w\|_2) dw
\]
which tends to 0 as $m$ goes to infinity.

As for $I\!I$, another integration by parts gives
\[
 I\!I= 
 2 \frac{\partial g}{\partial x}(1,1) \int_{\bR^{n-1}\setminus\{0\}} u(y) 
 h_m(\|y\|_2)  \frac{1}{\|y\|_2^2}
 \, dy 
\]
\[
 -2 \int_{\bR^{n-1}\setminus\{0\}} u(y) \int_0^{\|y\|_2} h_m(x) \left(
 \frac{\partial^2 g}{\partial x^2}(x,\|y\|_2) \frac{1}{g(x,\|y\|_2)^2} +
 \left(\frac{\partial g}{\partial x}(x,\|y\|_2) \right)^2 \frac{-2}{g(x,\|y\|_2)^3}
 \right)\, dx
  \, dy
\]
Here, $\frac{\partial g}{\partial x}(1,1)$ is defined as $\lim_{x \rightarrow 1^-} \frac{\partial g}{\partial x}(x,1)$. The first integral and the last term of the second integral are positive, so we only need to study the term
\[
  -2 \int_{\bR^{n-1}\setminus\{0\}} u(y) \int_0^{\|y\|_2} h_m(x) 
 \frac{\partial^2 g}{\partial x^2}(x,\|y\|_2) \frac{1}{g(x,\|y\|_2)^2}\, dx
  \, dy.
\]
Since $h_m(x)$ is an approximate identity and $\frac{\partial^2 g}{\partial x^2}(x,\|y\|_2)$ converges to zero uniformly in $y$ as $x$ goes to zero, this term tends to zero as $m$ goes to infinity and Lemma 1 holds.

\qed

\begin{remark} A similar proof shows that any centrally symmetric body (not necessarily of revolution) that has a cylindrical part is not an intersection body.
\end{remark}

\section{Regularity and convexity conditions for an intersection body of revolution to be the intersection body of a star body}

Let $L$ be a centered star body of revolution about the $x_n$-axis in $\mathbb{R}^n$, and let $\rho_L$ be its radial function, which we assume to be continuous. The radial function $\rho_L$ may be considered as a function of the angle $\varphi$ from the $x_n$-axis.

 Let $K$ be the intersection body of $L$,  defined as  the body whose radial function $\rho_K$ is the spherical Radon transform of $\rho_L^{n-1}/(n-1)$. Then 
\begin{equation}
  \label{radon1}
     \rho_K(\varphi)=\frac{2 \omega_{n-2}}{(n-1)\sin \varphi} \int_{\pi/2-\varphi}^{\pi/2}  \, \rho_L(\psi)^{n-1} \left(1-\frac{\cos^2 \psi}{\sin^2 \varphi} \right)^{(n-4)/2} \, \sin \psi \, d \psi,
\end{equation}
if $0<\varphi \leq \pi/2$, and $\rho_K(0)=\kappa_{n-1} \rho_L(\pi/2)^{(n-1)}$. Here, $\omega_{n}$  denotes the surface area of the unit ball in $\mathbb{R}^n$. 
 A derivation of this formula can be found in  \cite{G}, Theorem C.2.9. 

If we substitute $x=\sin \varphi$, $t=\cos \psi$ in (\ref{radon1}), we obtain 
\begin{equation}
 \label{radon2}
    \rho_K(\arcsin x)=\frac{2 \omega_{n-2}}{(n-1)x^{n-3}} \int_0^x \rho_L (\arccos t)^{n-1} (x^2-t^2)^{(n-4)/2}\, dt,
\end{equation}
for $0<x\leq 1$. When $\rho_L$ is continuous, this formula can be inverted:
\begin{equation}
  \label{invradon2}
    \rho_L (\arccos t)^{n-1}=\frac{1}{(n-3)! \omega_{n-1}} t \left( \frac{1}{t} \frac{d}{dt} \right)^{n-2} 
        \int_0^t  \rho_K(\arcsin x) x^{n-2} (t^2-x^2)^{(n-4)/2} \, dx,
\end{equation}
for $0 < t \leq 1$ (see \cite{G}, Corollary C.2.11).

Our first proposition determines  the relation between the regularity of $\rho_K$  and $\rho_L$, thus providing a necessary regularity condition for a body of revolution to be an intersection body of a star body.


\begin{Proposition} 
 \label{regularity}
 Let $L$ be a star body of revolution in $\mathbb{R}^n$,  where $n \geq 4$ is an even number. Assume that its radial function $\rho_L$  is of class $ C^m(S^{n-1})$. Let $K$ be the intersection body of $L$, with radial function $\rho_K(\varphi)$ given by (\ref{radon1}). Then $\rho_K(\varphi)$   is of class $C^{m+\frac{n}{2}-1}$  for $0<\varphi<\pi/2$, of class $C^m$ at $\varphi=0$, and of class $C^{m+n-2}$ at $\varphi=\pi/2$.
\end{Proposition}

\begin{Corollary}
 \label{necessary-reg}
Let $K$   be a body of revolution in $ \mathbb{R}^n$, with  $n \geq 4$ even.  A necessary condition for $K$ to be an intersection body of a star body is that its radial function $\rho_K(\varphi) $ is of class $C^{\frac{n}{2}-1}$ for $0<\varphi<\pi/2$, and of class $C^{n-2}$ at $\varphi=\pi/2$.
\end{Corollary}

Some immediate applications of Corollary \ref{necessary-reg} are the following:
\begin{itemize}
  \item Any body of revolution that is not $C^1$, such as the cylinder, or a cylinder with two conical caps, is not an intersection body of a star body in dimensions 4 and higher.
 \item A double cone is not an intersection body of a star body in dimensions 4 and higher, because its radial function is not $C^2$ at $\varphi=\pi/2$.

\end{itemize}

\bigskip

\noindent{\bf Proof of Proposition \ref{regularity}:}

{\bf Part 1:} We consider first the case $0<\varphi<\pi/2$. It will be more convenient to use equation (\ref{radon2}), with $0<x<1$.  Let us denote $r(t)= \rho_L (\arccos t)^{n-1}$, and  
\begin{equation}
  \label{rad}
  F(x)=x^{-n+3} \int_0^x r(t) (x^2-t^2)^{(n-4)/2}dt.
\end{equation}
 Thus, $ \rho_K(\arcsin x)= \frac{2 \omega_{n-2}}{(n-1)}F(x)$, and we have to prove that if $r(t) \in C^m$, then \newline $F(x) \in C^{m+\frac{n}{2}-1}$. 

The first derivative of $F$ is 
\[
       (-n+3)x^{-n+2} \int_0^x r(t) (x^2-t^2)^{(n-4)/2}dt + x^{-n+4}(n-4) \int_0^x r(t) (x^2-t^2)^{(n-6)/2}dt.
\]
Looking at the second term, we see that the exponent of $(x^2-t^2)$ has decreased by one. Continuing this process, the $k$-th derivative of $F$ will thus contain a term in which the exponent of $(x^2-t^2) $ is $(n-4-2k)/2$. Hence, the $(n-2)/2$-th derivative of $F$ will be the first to contain a term without integral, which is the term with the lowest regularity. It is equal to 
\[
     (n-4)!!\,  x^{-n/2+1} \, r(x).
\]
Thus, if $r \in C^m$, then $F \in C^{m+(n-2)/2}$.

\bigskip

{\bf Part 2:} To study the regularity at the point $x=0$, we extend $F$ evenly (since it is the radial function of a body of revolution). At $x=0$, the value of $F$ must be 
\[
   \sum_{j=0}^{(n-4)/2} \binom{\frac{n-4}{2}}{ j}  (-1)^j \frac{r(0)}{2j+1},
\]
so that $F$ is continuous. This is easy to see by expanding the term $(x^2-t^2)^{(n-4)/2}$ in Equation (\ref{rad}) and applying L'H\^opital's rule. 

We will show that, for every natural number $k$,
\begin{equation}
  \label{at0}
      \tilde{ F}^{(k)} (0+)=  \sum_{j=0}^{(n-4)/2} \binom{\frac{n-4}{2}}{ j}  (-1)^j  \lim_{x \rightarrow 0+} \frac{r^{(k)}(x)}{2j+k+1}.
\end{equation}
Hence the regularity of $F(x)$ at $x=0$ is the same as the regularity of $r(x)$ at $x=0$. 

We proceed by induction. Assume that (\ref{at0}) holds for every $k \leq k_0$. We will show that the formula is true for $k_0+1$. 
Expanding the binomial inside equation (\ref{rad}), we can write $F$ as 
\[
           F(x)=  \sum_{j=0}^{(n-4)/2} \binom{\frac{n-4}{2}}{ j}  (-1)^j    g_j(x) I_j(x),   
\]
where
\[
  g_j(x)=\frac{1}{x^{2j+1}} \,\, {\mbox{ and }}
     I_j(x)=\int_0^x r(t)t^{2j} \, dt \equiv \int_0^x r_j(t) \, dt 
\]
Then,
\begin{equation}
  \label{derivadas}
          F^{(k)}(x)=  \sum_{j=0}^{(n-4)/2} \binom{\frac{n-4}{2}}{ j}  (-1)^j  
             \left[ \sum_{i=0}^k \binom{k}{i}   g_j^{(i)}(x) I_j^{(k-i)}(x)
            \right],
\end{equation}
and 
\[
  F^{(k_0+1)}(0+)=\lim_{x \rightarrow 0+}  \frac{ F^{(k_0)}(x)-F^{(k_0)}(0)}{x}=
\]
\[
  = \lim_{x \rightarrow 0+} \sum_{j=0}^{(n-4)/2} \binom{\frac{n-4}{2}}{ j}  (-1)^j  
         \frac{1}{x}  \left[ \sum_{i=0}^{k_0} \binom{k_0}{i}   g_j^{(i)}(x) I_j^{(k_0-i)}(x)-\frac{r^{(k_0)}(0+)}{2j+k_0+1}      
 \right].
\]
Let us fix $j$ and consider the term 
\[
   T_j(x) =\frac{1}{x}  
             \left[ \sum_{i=0}^{k_0} \binom{k_0}{i}   g_j^{(i)}(x)  I_j^{(k_0-i)}(x)-\frac{r^{(k_0)}(0+)}{2j+k_0+1}
            \right]        .
\]
We have to prove that 
\begin{equation}
  \label{finish}
    \lim_{x \rightarrow 0} T_j(x) = \frac{r^{(k_0+1)}(0+)}{(2j+k_0+2)}.
\end{equation}
Observing that 
\[
      g_j^{(i)}(x)=\frac{(2j+i)!}{(2j)!} (-1)^i x^{-(2j+i+1)},
 \]
and, if $k_0-i \geq 1$,
\[
     I_j^{(k_0-i)}(x)=r_j^{(k_0-i-1)}(x) ,
\]
and multiplying both the numerator and denominator of $T_j(x)$ by $x^{2j+k_0+1}$, $T_j(x)$  can be rewritten as
\[
       T_j (x)=\frac{1}{x^{2j+k_0+2}}  
             \left[ \sum_{i=0}^{k_0-1} \binom{k_0}{i}  \frac{(2j+i)!}{(2j)!} (-1)^i x^{k_0-i} r_j^{(k_0-i-1)}(x)+ \right.
\]
\[
  \left. \frac{(2j+k_0)!}{(2j)!} (-1)^{k_0} \int_0^x r_j(t)\,dt-\frac{r^{(k_0)}(0+)}{2j+k_0+1}x^{2j+k_0+1}
            \right] .   
\]
By L'H\^opital's rule,
\[
     \lim_{x \rightarrow 0} T_j(x)=\lim_{x \rightarrow 0} \frac{1}{(2j+k_0+2)x^{2j+k_0+1}}  
             \left[ \sum_{i=0}^{k_0-1} \binom{k_0}{i}  \frac{(2j+i)!}{(2j)!} (-1)^i \left( (k_0-i) x^{k_0-i-1} r_j^{(k_0-i-1)}(x)+  \right. \right.
\]
\begin{equation}
 \label{limitTj}
  \left.  \left. x^{k_0-i} r_j^{(k_0-i)}(x) \right)+\frac{(2j+k_0)!}{(2j)!} (-1)^k_0 r_j(x)-r^{(k_0)}(0+)x^{2j+k_0}
            \right] 
\end{equation}
In the term inside square brackets, we group together the terms with derivatives  of $r_j(x)$ of the same order, obtaining
\[
    \left( 
         \sum_{i=0}^{k_0}  \frac{(2j+i-1)!}{(2j-1)!} (-1)^i r_j^{(k_0-i)}(x) x^{k_0-i}
    \right)-r^{(k_0)}(0+)x^{2j+k_0} .
\]
Recalling now that $r_j(x)=r(x)x^{2j}$, this sum equals
\begin{equation}
 \label{delta}
    2j \sum_{l=0}^{k_0-1} \binom{k_0}{l} r^{(l)}(x) x^{2j+l}
      \left( 
                \sum_{i=0}^{k_0-l} (-1)^i \binom{k_0-l}{i}\frac{(2j+i-1)!}{(2j+i-k_0+l)! }   
     \right) 
\end{equation}
\[ 
    + (r^{(k_0)}(x) -r^{(k_0)}(0+))x^{2j+k_0} .
\]
The inner sum appearing in (\ref{delta}) is zero for every value of $k_0$ and $l$, because it is equal to $\Delta^m P(x)$ evaluated at $x=2j$, where $P(x)=(x+1)(x+2)\cdots (x+m-1)$ and $\Delta$ is the difference operator $\Delta P(x)=P(x-1)-P(x)$.  But if we apply $\Delta^m$  to a polynomial of degree $m-1$, we obtain zero. Thus, $(\ref{delta})$ equals 
\[
   (r^{(k_0)}(x) -r^{(k_0)}(0+))x^{2j+k_0}.
\]
Writing this instead of the square brackets in  $(\ref{limitTj})$, we now have 
\[
      \lim_{x \rightarrow 0} T_j(x)=\frac{1}{(2j+k_0+2)}\, \lim_{x \rightarrow 0} \frac{(r^{(k_0)}(x) -r^{(k_0)}(0+))}{x}  
\]
\[
      =\frac{r^{(k_0+1)}(0+)}{(2j+k_0+2)}.
\]
This proves $(\ref{finish})$ and hence the regularity of $F$ at $0$ is the same as the regularity of $r$ at $0$.

\bigskip
 
{\bf Part 3:} Finally, we will show that $\rho_K(\varphi) \in C^{m+n-2}$  at $\varphi=\pi/2$. If we set $u=\pi/2-\psi$ and 
$s(u)=\rho_L(\pi/2-u)^{n-1} \left(1-\frac{\sin^2 u}{\sin^2 \varphi} \right)^{(n-4)/2} $, equation   $(\ref{radon1})$ becomes, disregarding the constants,
\[
     \rho_K(\varphi)=\csc  \varphi \, \int_0^\varphi s(u) \cos u \, du,.
\]
for $0 \leq \varphi \leq \pi/2$. Since $K$ is centrally symmetric, for $\pi/2<\varphi\leq \pi$  its radial function is  $\rho_K(\pi-\varphi)$. In particular, all even derivatives of $\rho_K$ are continuous at $\pi/2$, and the odd derivatives are continuous if and only if their value at $\pi/2$ is $0$. The same is true for the body $L$ and thus for $s(u)$.  Our hypothesis says that $s\in C^m$, and we will assume that  $s^{(m+1)}$ is not continuous at $\pi/2$.

Let  $I(\varphi)=\int_0^\varphi s(u) \cos u \, du$. Observe that $I'(\pi/2)=0$, and that all the odd derivatives of $g(\varphi)=\csc \varphi$ are 0 at $\varphi=\pi/2$. Hence,  $\rho_K'$ is continuous at $\pi/2$, and for $k \geq 3$ odd,
\[
     \rho_K^{(k)}\left(\frac{\pi}{2}\right)=\sum_{i=0}^{(k-3)/2} \binom{k}{2i}  g^{(2i)}\left(\frac{\pi}{2}\right) I^{(k-2i)}\left(\frac{\pi}{2}\right).
\]
The highest order derivative of $s(u)$  will appear in the term $I^{(k)}$. Since $s \in C^m \setminus C^{m+1}$ at $\varphi=\pi/2$, we have to find the first value of $k$ for which $\rho_K^{(k)}$ contains a term with  $s^{(m)}$  multiplied by a function which is non-zero at $\pi/2$. That value of $k$ will give us the regularity of $\rho_K$.

Since $s(u)=\rho_L(\pi/2-u)^{n-1} \left(1-\frac{\sin^2 u}{\sin^2 \varphi} \right)^{(n-4)/2} $, we need to differentiate the integral $I(\varphi)$  $(n-2)/2$ times in order for $s$ to appear outside of an integral. Indeed,  $ I^{((n-2)/2)}$ contains a term of the form  
\[
         \frac{ (\cos \varphi)^{(n-2)/2}}{(\sin \varphi)^{(n-4)/2}} \, s(\varphi).
\]
Let $ f(\varphi)= \frac{ (\cos \varphi)^{(n-2)/2}}{(\sin \varphi)^{(n-4)/2}}$. Then, $f^{(l)}(\pi/2)=0$ for $0 \leq l \leq (n-4)/2$, and $f^{((n-2)/2)}(\pi/2) \neq 0$. Thus, $ I^{(n-2+m)}$ is the first of the derivatives of $I$ containing the term $s^{(m)}(\varphi)f^{((n-2)/2)}(\varphi) $. This shows that $\rho_K$ has regularity $C^{m+n-2} \setminus C^{m+n-1}$ at $\pi/2$. 

\qed

\bigskip

Theorem 8.1.13 in \cite{G} proves that if a body $K \subset \mathbb{R}^4$ is axis-convex and its radial function $\rho_K(\varphi)$ is $C^1$ for $0< \varphi<\pi/2$ and $C^2$ at $\varphi=\pi/2$, then $K$ is  an intersection body of a star body. (The statement  of theorem actually asks for $\rho_K $ to be in $C^2$, but in its proof the continuity of the second derivative is only used at the point $\pi/2$). The regularity hypothesis ensures the continuity of the  inverse Radon transform of $\rho_K$, as in Proposition \ref{regularity}.  The  axis-convexity guarantees that the inverse Radon transform of $\rho_K$ is non-negative. It is well known that, in dimensions 5 and higher, there exist infinitely smooth convex bodies that are not intersection bodies of star bodies \cite{G}. Hence, we cannot expect a result similar to Theorem 8.1.13 in dimension $n \geq 5$. However, the following theorem shows that if a  body of revolution $K \subset \mathbb{R}^{2n_0}$ is an intersection body of a star body that verifies an additional axis-convexity-type property,  then $K$ is also an intersection body of a star body in dimension $2n_0+2$. 

Before stating the theorem, we will introduce some notation. We say that a body of revolution $L \subset \mathbb{R}^{n}$ is {\it equator-convex} if its intersection with every plane parallel to the equator of $L$ is convex. If we consider $L$ to be 2-dimensional,  $L$ is equator-convex if every line parallel to the x-axis intersects the body in a line segment.

Given a radial function $\rho_K(\varphi)$, $0\leq \varphi \leq \pi/2$, we will denote by $K_n$ the body of revolution in $\mathbb{R}^n$ whose radial function is $\rho_K$. 

\begin{Theorem}
 \label{equatorconvex}
Let $ n\geq 4$ even, and let  $K_n$  be a body of revolution with radial function $\rho_K$. Assume that $K_n$ is the intersection body of a star body $L \subset \mathbb{R}^{n}$. If  $L$ is equator-convex and $\rho_L \in C^1$ , then $K_{n+2}$ is also an intersection body of a star body. If $L$ is not equator-convex, then $K_{n+2}$ is not an intersection body. 
\end{Theorem}

{\bf Proof:} Since $K_n$ is the intersection body of $L$, and denoting $f(t)= \rho_L (\arccos t)$, we have  by (\ref{radon2}): 
\begin{equation}
 \label{BIL}
     \rho_K(x)=\frac{2 \omega_{n-2}}{(n-1)x^{n-3}} \int_0^x f(t)^{n-1} (x^2-t^2)^{(n-4)/2}\, dt.
\end{equation}
Let $\widetilde{\rho}(t)$ be the  $(n+2)$-inverse Radon transform of  $\rho_K$. To show that $K_{n+2}$ is an intersection body of a star body, we need to check that $\widetilde{\rho}(t)$  is a non-negative continuous function. By (\ref{invradon2}) and (\ref{BIL}),
\[
    \widetilde{\rho}(t)= \frac{1}{(n-1)! \, \omega_{n+1}} t \left( \frac{1}{t} \frac{d}{dt} \right)^{n} 
        \int_0^t  \rho_K(\arcsin x) x^{n} (t^2-x^2)^{(n-2)/2} \, dx
\]
\begin{equation}
 \label{n+2}
      = c_n \, t \left( \frac{1}{t} \frac{d}{dt} \right)^{n} 
        \int_0^t  x^{3} (t^2-x^2)^{(n-2)/2}  
         \left( \int_0^x f(u)^{n-1} (x^2-u^2)^{(n-4)/2}\, du \right)  \, dx,
\end{equation}
where $c_n=\frac{2\omega_{n-2}}{(n-1)! \, \omega_{n+1} (n-1)}$.
After $(\frac{n}{2})$ applications of the operator $\left( \frac{1}{t} \frac{d}{dt} \right)$, equation (\ref{n+2}) becomes
\begin{equation}
 \label{aftern/2}
      (n-2)!!\, c_n \, t \left( \frac{1}{t} \frac{d}{dt} \right)^{\frac{n}{2}}  t^2
        \int_0^t f(u)^{n-1} (t^2-u^2)^{(n-4)/2}\, du .
\end{equation}
The next application of $\left( \frac{1}{t} \frac{d}{dt} \right)$ results in two terms,
\[
      (n-2)!!\, c_n \, t \left( \frac{1}{t} \frac{d}{dt} \right)^{\frac{n}{2}-1}  \left[ 
       2 \int_0^t f(u)^{n-1} (t^2-u^2)^{(n-4)/2}\, du  \right.
\]
\[
  \left. +
 (n-4) \, t^2 \int_0^t f(u)^{n-1} (t^2-u^2)^{(n-6)/2}\, du
    \right].
\]
Let us call 
\[
    A=(n-2)!!\, c_n \, t \left( \frac{1}{t} \frac{d}{dt} \right)^{\frac{n}{2}-1} 
       2 \int_0^t f(u)^{n-1} (t^2-u^2)^{(n-4)/2}\, du  
\]
and 
\[
     B=(n-2)!!\, c_n \, t \left( \frac{1}{t} \frac{d}{dt} \right)^{\frac{n}{2}-1}   (n-4) \,t^2 \int_0^t f(u)^{n-1} (t^2-u^2)^{(n-6)/2}\, du
\]
It is not hard to see that 
\[
     A=2 c_n (n-2)!! (n-4)!! f(t)^{n-1}.
\]
As for $B$, it is the same as  equation (\ref{aftern/2}), with $n$ replaced by $n-2$. Hence, after $(\frac{n}{2}-3)$ more differentiations, we obtain
\[
    A+B=c_n (n-2)!! (n-4)!!\left[ (n-4) f(t)^{n-1}+t \left(\frac{d}{dt}\right)^2 t^2 \int_0^t f(u)^{n-1} \, du \right]
\]
\[ 
    =c_n (n-2)!! (n-4)!!\left[ (n-4) f(t)^{n-1}+ \frac{d}{dt} \left( 2\int_0^t f(u)^{n-1} \, du + t f(t)^{n-1} \right) \right]
\]
\[ 
    =c_n (n-2)!! (n-4)!!\left[ (n-1) f(t)^{n-1}+ (n-1)tf(t)^{n-2} f'(t) \right]
\]
\[ 
    =\frac{1}{2\pi}  f(t)^{n-1} \left[ 1 + \frac{t f'(t)}{f(t)} \right].
\]
Thus,  $\widetilde{\rho}(t)= \frac{1}{2\pi}  f(t)^{n-2} \left[ f(t) + t f'(t) \right]$. Since $f \in C^1$, $\widetilde{\rho}$ is continuous. Since $L$ is equator-convex, $(t f(t))$ is increasing, which means that $f(t) + t f'(t)$ and  $\widetilde{\rho}$ are non-negative. We conclude that $K_{n+2}$ is an intersection body of the body whose radial function is $\left( (n+1) \widetilde{\rho}(t)  \right)^{1/(n+1)}$.

On the other hand, if $L$ is not equator-convex, then $\widetilde{\rho}$ takes negative values and $K_{n+2}$ is not an intersection body (with the general definition). In this case, the assumption $f  \in C^1$ may be relaxed.

\qed

\begin{ex}
  \label{diabolo}
\end{ex}
Our first application of Theorem  \ref{equatorconvex} will be the construction of a body that is an intersection body of a star body up to dimension $2n_0$, and is not an intersection body starting from dimension $2n_0+2$. Consider the star body of revolution  $L$ whose radial function is 
\[
     \rho_L(\psi)=\left\{  
        \begin{array}{ll}
           \displaystyle    \frac{2\cos \psi + \sin \psi}{5\cos^2 \psi -1} , &\,\, 0\leq \psi<\pi/4 \\
              \csc \psi, & \,\, \pi/4<\psi \leq \pi/2 \\
        \end{array}.
    \right.
\]

\begin{figure}[!b]
 \begin{center}
  \includegraphics[width=1.5in]{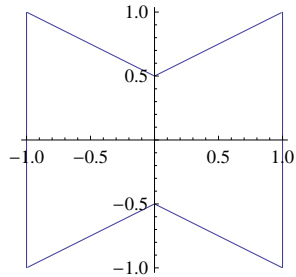}
 \end{center}

\caption{\small Cross-section of the body $L$ in Example $\ref{diabolo}$.}
\label{fig-label1}
\end{figure}

\begin{figure}[!b]
\begin{center}
\includegraphics[width=1.5in]{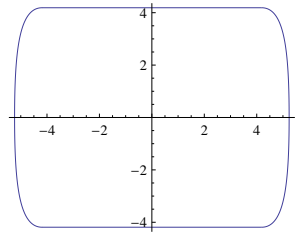}
\end{center}

\caption{\small Cross-section of the body $K_4=I \! L$ in Example $\ref{diabolo}$.  In dimension 6, this body is not an intersection body.}
\label{fig-label2}
\end{figure}

\begin{figure}[!b]
  \begin{center}
    \includegraphics[width=2in]{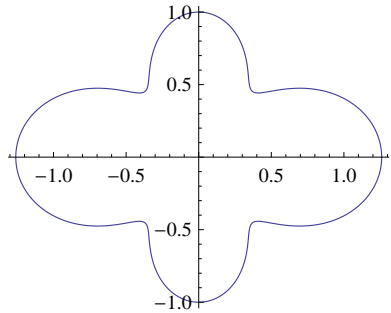}
  \end{center}

  \caption{\small Cross-section of the body $\widetilde{L}$ in Example \ref{diabolo}. }
  \label{fig-label3}
\end{figure}


The function $\rho_L$ is just continuous, and $L$ is not equator-convex (see Figure \ref{fig-label1}). Let 
$K_{2n_0}$ be the intersection body of $L$ in dimension $n=2n_0$, and $\rho_K$ its radial function.  Notice that, although $\rho_L$ is not $C^1$, its derivative is piecewise continuous. Hence, in dimension $2n_0+2$, $\rho_K$  is the Radon transform of a piecewise continuous, sign-changing function., and this means that $K_{2n_0+2}$ is not an intersection body.  Figure \ref{fig-label2} shows the cross-section of the intersection body of $L$ in dimension 4.

We can also start with a body  that has higher regularity. For example, let $\rho_{\widetilde{L}}(\psi)=\left(2 - 6 \cos^2 \psi +5 \cos^4 \psi \right )^{1/3}$. This function is $C^\infty$, but $\widetilde{L}$ is not equator-convex (see Figure \ref{fig-label3}). As before,  let $\rho_{\widetilde{K}}$ be given by (\ref{radon1}) with $\rho_L=\rho_{\widetilde{L}}$ and $n=2n_0$ . Then, by Theorem  \ref{equatorconvex}, $\widetilde{K}_{2n_0}$  is an intersection body of a star body, but $\widetilde{K}_{2n_0+2}$ is not an intersection body. 

\begin{figure}[!b]
  \begin{center}
    \includegraphics[width=2in]{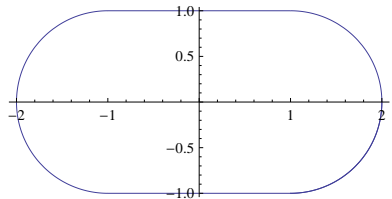}
  \end{center}

  \caption{\small Cross-section of Lonke's barrel zonoid (see Example \ref{lonke}). }
  \label{fig-label4}
\end{figure}

\begin{ex}
  \label{lonke}
\end{ex}
We will now study the barrel $B=B_n + B_{n-1} \subset \mathbb{R}^n$, where $B_n$ is the unit ball in $\mathbb{R}^n$. See Figure \ref{fig-label4}. This body was introduced by Lonke in the paper \cite{L}, where he proved that, in dimensions 3 and 4, $B$ is a zonoid whose dual is a zonoid (and, in particular, $B$ is an intersection body). Lonke's proof that the dual of $B$ is a zonoid does not work in dimensions 5 and higher, and we are interested in studying if the reason can be that $B$ no longer is an intersection body. 

\begin{figure}[!b]
  \begin{center}
    \includegraphics[width=1in]{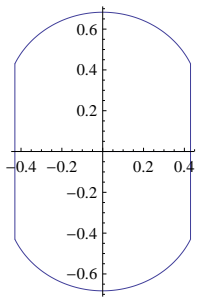}
  \end{center}

  \caption{\small Cross-section of the body whose intersection body in $\bR^4$ is Lonke's barrel zonoid (see Example \ref{lonke}). }
  \label{fig-label5}
\end{figure}

 Its  radial function is 
\[
     \rho_B(\varphi)=\left\{  
        \begin{array}{ll}
               \sec\varphi, &\,\, 0\leq \varphi<\pi/4 \\
              2\sin \varphi, & \,\, \pi/4<\varphi \leq \pi/2 \\
        \end{array},
    \right.
\]
and thus $\rho_B(\varphi)$ is $ C^1$  at the point $\varphi=\pi/4$, and  $C^\infty$ everywhere else. By Proposition \ref{regularity}, $B$ is not an intersection body of a star body in dimensions 6 and higher. In dimension 4, we use  the inversion  formula (\ref{invradon2}) to obtain that $B$ is the intersection body of the body $L$ whose radial function is
\[
     \rho_L(\varphi)=\left\{  
        \begin{array}{ll}
           \displaystyle    \left( \frac{3\cos\psi}{\pi} \right)^{1/3},  &\,\, 0\leq \psi<\pi/4 \\ \\
            \displaystyle  \left(\frac{3}{4\pi} \right)^{1/3} \csc \psi,  & \,\, \pi/4<\psi \leq \pi/2 \\
        \end{array},
    \right.
\]
Note that $\rho_L'$ is piecewise continuous and that the body $L$ is equator-convex (see Figure \ref{fig-label5}). Although $B$ is not an intersection body of a star body in dimension 6,  it is an intersection body in $\mathbb{R}^6$, since $\rho_B$ is the Radon transform of a non-negative piecewise continuous function. 

Proposition \ref{regularity} shows that, for a given function $\rho_K$, increasing the dimension by two units decreases the regularity of its inverse Radon transform by 1. Since the inverse Radon transform of $\rho_B(\varphi)$ is piecewise continuous in dimension 6, we expect that in dimension 8 it will contain delta functions. Indeed,   in the sense of distributions  $\rho_B(\arcsin x)$ equals, up to a constant, the Radon transform of $\rho_{L_8}(\arccos t)-\frac{4}{15} \delta_{1/\sqrt{2}}(t)$, where
$\delta_{1/\sqrt{2}}$ is the Dirac measure supported at the point $t=1/\sqrt{2}$, and
\[
    \rho_{L_8}(\arccos t)=\left\{  
        \begin{array}{ll}
             \displaystyle   \frac{1}{(1-t^2)^{7/2}},  & \,\, 0< t < 1/\sqrt{2} \\
            \frac{96}{15} t,  &\,\,  1/\sqrt{2} \leq t < 1 \\ 
        \end{array}.
    \right.
\]
Since this measure is negative at the point $t=1/\sqrt{2}$, $B$ is not an intersection body in dimensions 8 and higher. Thus, the dual of Lonke's barrel's is not a zonoid in dimensions 8 and higher. Only in dimensions 5, 6 and 7 the question is still unanswered.


\begin{thebibliography}{99} 

\bibitem{B} E. D. Bolker,
\newblock \emph{A class of convex bodies,}
\newblock Trans. Amer. Math. Soc. 145 (1969), 323--345.


\bibitem{G} R. J. Gardner,
\newblock \emph{Geometric Tomography, 2nd edition}, 
\newblock Cambridge University Press, 2006.


\bibitem{GKS} R. J. Gardner, A. Koldobsky, T. Schlumprecht,
\newblock \emph{An analytic solution to the Busemann-Petty problem on sections of convex bodies,}
\newblock Ann. Math. 149 (1999), 691--703.

\bibitem{GW} P. Goodey, W. Weil,
\newblock \emph{Intersection bodies and ellipsoids}, 
\newblock Mathematika 42 (1995), 295--304.


\bibitem{K} A. Koldobsky,
\newblock \emph{Fourier Analysis in Convex Geometry,}
\newblock Math. Surveys and Monographs, AMS (2005).





\bibitem{K5} A. Koldobsky,
\newblock \emph{Second derivative test for intersection bodies,}
\newblock Adv. Math. 136 (1998), 15--25.





\bibitem{LP} J. Lindenstrauss, A. Pelczy\'nski,
\newblock \emph{Absolutely summing operators in $L_p$ spaces and their applications,}
\newblock Studia Math. 29 (1968), 257--326.

\bibitem{L} Y. Lonke,
\newblock \emph{On zonoids whose polars are zonoids,}
\newblock Israel J. Math. 102 (1997), 1--12.

\bibitem{Lu} E. Lutwak,
\newblock \emph{Intersection bodies and dual mixed volumes,}
\newblock Advances in Math. 71 (1988), 232--261.

\bibitem{S} R. Schneider,
\newblock \emph{Zonoids whose polars are zonoids,}
\newblock Proc. AMS 50 (1974), 365--368.

\bibitem{Sc} R. Schneider, 
\newblock \emph{Convex bodies: The Brunn-Minkowski Theory,}
\newblock Cambridge University Press, 1993.

\bibitem{Z} G. Zhang,
\newblock \emph{Intersection bodies and polytopes,}
\newblock Mathematika, 46 (1999), 29--34.

\end{thebibliography}
\end{document}